\documentclass[MSNbibl,number,citesort,seceqn,dvips]{arxbj}
\usepackage{dcolumn}

\aid{0}
\volume{19}
\issue{1}
\pubyear{2013}
\firstpage{228}
\lastpage{251}
\doi{10.3150/11-BEJ400}

\makeatletter
\def\bsuffix #1{#1}
\newcolumntype{d}[1]{D{.}{.}{#1}}
\newtheorem{lemma}{Lemma}
\newtheorem{theorem}{Theorem}
\makeatother

\begin{document}
\begin{frontmatter}

\title{Parameter estimation and model testing for Markov processes via
conditional characteristic functions}
\runtitle{CCF-based estimation and testing}

\begin{aug}
\author[1,2]{\fnms{Song X.} \snm{Chen}\thanksref{1,2,e1}\ead[label=e1,mark,text=songchen@iastate.edu]{songchen@iastate.edu}},
\author[3]{\fnms{Liang} \snm{Peng}\corref{}\thanksref{3}\ead[label=e2,text=peng@math.gatech.edu]{peng@math.gatech.edu}}
\and
\author[1]{\fnms{Cindy L.} \snm{Yu}\thanksref{1,e3}\ead[label=e3,mark]{cindyyu@iastate.edu}}
\runauthor{S.X. Chen, L. Peng and C.L. Yu} 
\address[1]{Department of Statistics, Iowa State University, Ames, IA
50011-1210, USA.\\ \printead{e1,e3}}
\address[2]{Guanghua School of Management and Center for Statistical
Science, Peking University, Beijing 100871, China}
\address[3]{School of Mathematics, Georgia Institute of Technology,
Atlanta, GA 30332, USA.\\ \printead{e2}}
\end{aug}

\received{\smonth{5} \syear{2011}}
\revised{\smonth{8} \syear{2011}}

%
\begin{abstract}
Markov processes are used in a wide range of disciplines, including finance.
The transition densities of these processes are often unknown.
However, the conditional characteristic functions are more likely to be
available, especially for L\'evy-driven processes.
We propose an empirical likelihood approach, for both parameter
estimation and model specification testing, based on the conditional
characteristic function for processes with either continuous or
discontinuous sample paths. Theoretical properties of the empirical
likelihood estimator for parameters and a smoothed empirical likelihood
ratio test for a parametric specification of the process are provided.
Simulations and empirical case studies are carried out to confirm the
effectiveness of the proposed estimator and test.
\end{abstract}

%

\end{frontmatter}

\section{Introduction}\label{sec1}

Let $\{X_t(\theta)\}_{t \in{\mathcal{T}}}$ be a parametric $d$-dimensional
Markov process
defined by
\begin{equation}
\mathrm{d}X_t=\mu(X_t; \theta)\,\mathrm{d}t + \sigma(X_t; \theta) \,\mathrm{d} L_{t; \theta},
\label{eq:paralevy}
\end{equation}
where $\mu(\cdot)$ is a $d$-dimensional drift function, $\sigma
(\cdot)$ is a $d \times d$ matrix-valued function of $X_t$,
$L_{t;\theta
}$ is a L\'evy process in $R^d$
and $\theta\in\Theta\subset R^p$. When $L_t$ is a standard Brownian
motion, (\ref{eq:paralevy}) is a diffusion process having a continuous
sample path. When $L_t$ contains the Brownian motion and a compound
Poisson process, (\ref{eq:paralevy}) becomes the jump diffusion process.
A stochastic process of form (\ref{eq:paralevy}) has long been used to
model stochastic systems arising in physics, biology and other natural
sciences. It has also been the fundamental tool in financial modeling.
We refer to Sundaresan~\cite{Sun00} and Fan~\cite{Fan05} for overviews,
Barndorff-Nielsen, Mikosch and Resnick~\cite{Ba01} for recent developments
on L\'evy-driven processes and S{\o}rensen~\cite{Sre91} for statistical inference.
Important subclasses of (\ref{eq:paralevy}) include (i)
the multivariate diffusion process defined by
\begin{equation}
\mathrm{d}X_t=\mu(X_t; \theta)\,\mathrm{d}t + \sigma(X_t; \theta) \,\mathrm{d} B_t,
\label{eq:1.1}
\end{equation}
where $B_t$ is the standard Brownian motion in $R^d$ (Stroock and
Varadhan~\cite{StrVar79} and {\O}ksendal~\cite{kse00}); (ii) the Vasicek with Merton Jump
model (VSK-MJ) defined by
\begin{equation}\label{eq:VSK-MJ}
\mathrm{d}X_{t} = \kappa(\alpha-X_{t})\,\mathrm{d}t + \sigma \,\mathrm{d}B_{t}+J_{t}\,\mathrm{d}N_{t},
\end{equation}
where $\kappa$, $\alpha$ and $\sigma$ are unknown parameters and
represent the mean reverting rate, long-run mean and volatility of the
process, respectively, $N_{t}$ is a Poisson process with intensity
$\lambda$ and $J_{t}$ is the random jump size independent of the
filtration ${\mathcal F}_{t}$ up to time $t$ and has a normal density $N(0,
\eta^2)$ (Merton~\cite{Mer76}); (iii)
L\'evy driven Ornstein--Uhlenbeck process defined by
\begin{equation}\label{eq:IG-OU}\mathrm{d}X_{t} = -\lambda X_{t} \,\mathrm{d}t +
\mathrm{d}L_{\lambda t}, \qquad X_{0}>0,
\end{equation}
where $L_{t}$ is a L\'evy process with no Brownian part, a non-negative
drift and a L\'evy measure which is zero on the negative half line, and
the parameter $\lambda$ is positive (see Barndorff-Nielsen and
Shephard~\cite{BarShe01}).

Often a closed form expression for the transition density of process
(\ref{eq:paralevy}) is not available except for
some special processes, even if the transition density exists and is unique.
This fact prevents the use of the maximum
likelihood estimation (MLE) and the specification tests based on the
exact transition density.
Recently A\"{\i}t-Sahalia~\cite{AtS02,AtS08} established expansions for the transition
densities so that parameter estimation could be based on the
approximate likelihood functions. Testing may be also formulated via
the approximate density; see Chen, Gao and Tang~\cite{CheGaoTan08} and A\"{\i
}t-Sahalia, Fan
and Peng~\cite{AtSFanPen09} for such tests.
The conditional characteristic functions (CCF) are more likely
available than the transition densities for the continuous-time models,
especially for the L\'evy-driven processes through the celebrated L\'
evy--Khintchine representation. For instance, Duffie, Pan and Singleton
\cite{DufPanSin00} derived the explicit form of the CCF for multivariate affine
jump processes, which include the Vasicek with Merton jump process
given in (\ref{eq:VSK-MJ}). The CCF for the L\'evy-driven
Ornstein--Uhlenbeck process (\ref{eq:IG-OU}) is established in
Barndorff-Nielsen and Shephard~\cite{BarShe01}.

Statistical inference based on the characteristic functions was
proposed by Feuerverger and Mureika~\cite{FeuMur77},
Feuerverger and McDunnough~\cite{FeuMcD81} for independent observations and
Feuerverger~\cite{Feu90} for discrete time series.
Singleton~\cite{Sin01} introduced the approach to
inference for parametric continuous-time Markov processes and showed
that estimation can be carried out based on the CCF without having to
carry out the
the Fourier inversion. Chacko and Viceira~\cite{ChaVic03} proposed a generalized
method of moment estimator (GMM) for parameters at a finite number
of frequencies of the CCF. Carrasco \textit{et al.}~\cite{Caretal07} carried out GMM
estimation on a slowly diverging
number of frequencies of the CCF to achieve the optimal estimation
efficiency offered by the MLE. Jiang and Knight~\cite{JiaKni02} proposed GMM
estimators based on the joint characteristic function of the observed
state variables. Chen and Hong~\cite{CheHon10} proposed a test for multivariate
processes based on the CCF via a generalized spectral density approach.

In this paper, we first propose an empirical likelihood (Owen~\cite{Owe88})
approach for parameter estimation and model specification testing of a
parametric Markov process via the CCF. An empirical likelihood ratio is
formulated for the unknown parameters assuming specification (\ref
{eq:paralevy}), which leads to a non-parametric maximum likelihood
estimator. The proposed estimator may be viewed as a compromise between
Chacko and Viceira's~\cite{ChaVic03} GMM, based on a finite number of
frequencies, and that of Carrasco \textit{et al.}~\cite{Caretal07}, of a
high-dimensional GMM. The high-dimensional GMM approach requires
ridging a high-dimensional weighting matrix in order to avoid its
singularity, and the selecting the ridging parameter can be
computationally expensive. The proposed estimation utilizes a wide
range of frequency information in the parametric CCF, while having the
computation easily managed.

We then formulate an empirical likelihood CCF-based model specification
test for the parametric process (\ref{eq:paralevy}) via kernel
smoothing. The proposed test extends the transition density based tests
of Qin and Lawless~\cite{Qin94}, Chen, Gao and Tang~\cite{CheGaoTan08} and A\"{\i}t-Sahalia,
Fan and Peng
\cite{AtSFanPen09} to the CCF based. This largely increases the range of the
continuous-time Markov processes which can be tested directly without
replying on the transition density approximation. The proposed test
provides an alternative formulation of the CCF-based test of Chen and
Hong~\cite{CheHon10}, which is based on an explicit $L_2$ measure between an
kernel estimator of the CCF and its parametric counter-part. It is
largely distinct from the above mentioned tests, except Chen and Hong
\cite{CheHon10}, by targeting directly on CCF, which is more readily available
for continuous-time models than the transition density functions.
Another advantage of the proposed test is the empirical likelihood (EL)
formulation, which can produce an integrated likelihood ratio test in a
nonparametric setting. The proposed test utilizes some of the
attractive properties of the EL, like internal studentizing without an
explicit variance estimation and good power performance.
How to extend the proposed methods to the case of latent variables is
quite challenging and will be a part of our future research.

The paper is organized as follows. In Section~\ref{sec2}, we introduce and
evaluate the CCF-based empirical likelihood\ estimator. The model
specification test is
given in Section~\ref{sec3}. Section~\ref{sec4} reports results from simulation studies.
An empirical study for a set of 3-month treasury bill rate data is
analyzed in Section~\ref{sec5}.
All technical details are reported in the \hyperref[app]{Appendix}.

\section{Parameter estimation}\label{sec2}

Let $\{X_{t \delta}\}_{t=1}^n$ be $n$ discretely sampled observations
of (\ref{eq:paralevy}).
For notation simplification, we denote $X_{t \delta}$ as $X_t$,
where the sampling interval $\delta$ is any fixed quantity.
Let $\psi_t(u; \theta) = E_{\theta}(\mathrm{e}^{\mathrm{i} u^T X_{t+1} }| X_t)$, for $u
\in R^d$, be the conditional characteristic function. We use $\bar a$
and $A^{\star}$ to denote the conjugate of a complex number $a$ and the
conjugate transpose of the complex matrix $A$, respectively.\vadjust{\goodbreak}

Let $\varepsilon_t(\tau;\theta) = w(u,r; X_t) \lbrace \mathrm{e}^{\mathrm{i} u^T X_{t+1}} -
\psi_t(u; \theta)\rbrace$ for $\tau=(u^T, r^T)^T \in R^{2d}$, where
$w(u,r;X_t)$ is a weight factor. Here $\varepsilon_t(\tau;\theta)$ can be
regarded as ``residuals'' between $\mathrm{e}^{\mathrm{i} u^T X_{t+1}}$ and the parametric
CCF $\psi_t(u; \theta)$. The complex weight factor $w(u,r;X_t)$
satisfies $\bar{w}(u, r; X_t) = w(-u, -r; X_t)$ and $|w(u, r; X_t)| =
1$ for any $u, r \in R^d$, whose use is aimed to utilize more model information.
Let $\theta_0$ be the true parameter and the unique solution of
\begin{equation}
E\{ \mathrm{e}^{\mathrm{i} u^T X_{t+1}} - \psi_t(u; \theta) |X_t \}=0\qquad \mbox{for all
$u \in R^{d}$}. \label{eq:mean}
\end{equation}

From the Markov property and (\ref{eq:mean}), for any $\tau=(u^T,
r^T)^T \in R^{2d}$,
\begin{eqnarray}
&& E\{ \varepsilon_t(\tau; \theta_0)\} = 0\quad \mbox{and}\quad \operatorname{Cov}\{
\varepsilon_{t_1} (\tau; \theta_0), \varepsilon_{t_2} (\tau; \theta_0)
\} = 0
\qquad\mbox{if } t_1 \ne t_2. \label{eq:cov}
\end{eqnarray}

Let $\varepsilon_t^R(\tau; \theta)$\vspace*{1pt} and $\varepsilon_t^I(\tau;\theta)$
be the
real and imaginary parts of $\varepsilon_t(\tau; \theta)$ respectively,
and $\vec{\varepsilon }_t(\tau; \theta) = ( \varepsilon_t^R(\tau;
\theta),
\varepsilon_t^I(\tau; \theta) )^T$ be the real bivariate vector
corresponding to $\varepsilon_t(\tau; \theta)$.

We now formulate an empirical likelihood\ for $\theta$ based on the
CCF $\psi_t(u; \theta
)$. The empirical likelihood (EL) introduced in Owen~\cite{Owe88} is a
technique that allows construction of a non-parametric likelihood for
parameters of interest. Despite that the EL method is intrinsically
non-parametric, it possesses two important properties of a parametric
likelihood, the Wilks theorem and the Bartlett correction;
see Chen and Van Keilegom~\cite{CheVan09} for a latest overview and Kitamura,
Tripathi and Ahn~\cite{KitTriAhn04} for a formulation with conditional moments.

Let $p_1(\tau), \ldots, p_n(\tau)$ be probability weights allocated to
the ``residuals'' $\{\vec{\varepsilon }_t(\tau;\theta)\}_{t=1}^n$. A
local EL
for $\theta$ at $\tau$ is
\begin{equation}
L_n(\tau, \theta) = \max\prod_{t=1}^n p_t(\tau) \label{eq:EL},
\end{equation}
subject to $\sum_{t=1}^n p_t(\tau) = 1$ and $\sum_{t=1}^n p_t(\tau)
\vec
{\varepsilon }_t(\tau; \theta) = 0$.
Here the second constraint reflects (\ref{eq:mean}). The maximum
empirical likelihood\
is attained at $p_t(\tau) \equiv n^{-1}$ for all $t$ such that the
maximum likelihood $L_n(\tau;\theta) = n^{-n}$. Let $\ell_n(\tau;
\theta
) = - 2 \log\{L_n(\tau; \theta)/n^{-n}\}$ be the local log-EL ratio of
$\theta$ at $\tau$.

Employing the EL algorithm (Owen~\cite{Owe88}), the optimal $p_t(\tau)$ of the
above optimization problem (\ref{eq:EL}) is
\[
p_t(\tau) = {1 \over n} {1\over1 + \lambda(\tau;\theta)^T \vec
{\varepsilon
}_t(\tau; \theta) },
\]
where $\lambda(\tau;\theta)$ is a Lagrange multiplier in $R^2$ that satisfies
\begin{equation}
Q_{1 n}(\tau; \theta, \lambda) =: {1 \over n} \sum_{t=1}^n {\vec
{\varepsilon
}_t(\tau; \theta) \over1 + \lambda(\tau;\theta)^T \vec{\varepsilon
}_t(\tau;
\theta) } = 0. \label{eq:3.5.2}
\end{equation}
Hence, the local EL ratio becomes
\begin{equation}
\ell_n(\tau; \theta) = 2 \sum_{t=1}^n \log\lbrace1 + \lambda
(\tau
;\theta)^T \vec{\varepsilon }_t(\tau; \theta)\rbrace. \label{eq:3.5.2b}
\end{equation}
Integrating $\ell_n(\tau; \theta)$ against a probability weight $\pi
(\tau)$, which is supported on a compact set $S$ in $R^{2d}$, an
integrated empirical likelihood\ ratio for $\theta$ is
\begin{equation}
\ell_n(\theta) = \int_{\tau\in R^{2d} } \ell_n(\tau; \theta) \pi
(\tau)
\,\mathrm{d} \tau. \label{eq:intELratio}
\end{equation}
The maximum EL estimator (MELE) for $\theta$ is defined as
\[
\hat{\theta}_n = \arg\min_{\theta} \ell_n(\theta),
\]
by noting that $-2$ has been multiplied in the EL ratio $\ell_n(\tau;
\theta)$.

Like Qin and Lawless~\cite{Qin94}, we first show that there exists a
consistent estimator $\hat\theta_n$ with a certain rate of convergence
as follows.

\begin{lemma}
\label{le1}
Under Conditions \textup{C1--C4} given in the \hyperref[app]{Appendix}, with probability one,
$\ell_n(\theta)$ attains its minimum at $\hat\theta_n$ in the interior
of the ball $\Vert \theta-\theta_0\Vert \le \mathrm{O}(n^{-1/3})$, and $\hat\theta_n$
and $\lambda(\tau;\hat\theta_n)$ satisfy
\begin{equation}
\label{th1-1}
\cases{
 Q_{1n}(\tau;\hat\theta_n,\lambda(\tau;\hat\theta_n))=0 \qquad \mbox{for
all } \tau\in S \quad\mbox{and} \vspace*{2pt}\cr
\displaystyle\int Q_{2n}(\tau;\hat\theta_n,\lambda(\tau;\hat\theta_n))\pi
(\tau)
\,\mathrm{d}\tau=0,}
\end{equation}
where $Q_{1n}$ is defined in (\ref{eq:3.5.2}) and
\begin{equation}
\label{th1-2}
Q_{2n}(\tau;\theta,\lambda)=\frac1n\sum_{t=1}^n\frac1{1+\lambda
(\tau
;\theta)^T\vec{\varepsilon }_t(\tau;\theta)}\frac{\partial\vec
{\varepsilon }^T_t(\tau
;\theta)}{\partial\theta}\lambda.
\end{equation}
\end{lemma}

Before deriving the asymptotic normality of the $\hat\theta_n$, we define
\[
M_0 = {1\over2} \pmatrix{1 & 1  \vspace*{2pt}\cr i^{-1} & - i^{-1} }
, \qquad\tilde{\varepsilon }_t(\tau; \theta) = ( \varepsilon _t(\tau; \theta),
\varepsilon _t(-\tau; \theta) )^T,
\]
$A(\tau_1, \tau_2; \theta_0, \theta) = \operatorname{Cov}\{ \tilde{\varepsilon
}_1(\tau_1;
\theta), \tilde{\varepsilon }_1(\tau_2; \theta) \}$,
\begin{equation}
\Gamma(\theta_0) =: \int E\biggl ( { \partial\tilde{\varepsilon}_1^{\star}
(\tau; \theta_0) \over\partial\theta} \biggr) A^{-1} (\tau,\tau;
\theta_0,\theta_0) E \biggl( { \partial\tilde{\varepsilon}_1(\tau; \theta
_0) \over\partial\theta} \biggr) \pi(\tau) \,\mathrm{d} \tau\label{eq:Gamm1}
\end{equation}
and
\begin{eqnarray*}
V(\theta_0) &=& \int\int E \biggl( { \partial\tilde{\varepsilon}_1^{\star
} (\tau_1; \theta_0) \over\partial\theta} \biggr) A^{-1}(\tau_1,\tau
_1;\theta_0, \theta_0) A(\tau_1, \tau_2; \theta_0, \theta_0) \\
&&\phantom{\int\int }\times A^{\star-1}(\tau_2,\tau_2;\theta_0, \theta_0)E \biggl( {
\partial\tilde{\varepsilon}_1 (\tau_2; \theta_0) \over\partial
\theta}
\biggr) \pi(\tau_1) \pi(\tau_2) \,\mathrm{d} \tau_1 \,\mathrm{d} \tau_2.
\end{eqnarray*}

\begin{theorem}\label{th1} Under Conditions \textup{C1--C4} given in the
\hyperref[app]{Appendix}, for the estimator $\hat\theta_n$ in Lemma~\ref{le1}, we have
$\sqrt{n} (\hat{\theta}_{n} - \theta_0) \stackrel{d} \to N(0,
\Sigma)$
where $\Sigma= \Gamma^{-1}(\theta_0) V(\theta_0) \Gamma
^{-1}(\theta
_0)$.
\end{theorem}

The proposed estimator attains the $\sqrt{n}$-rate of convergence. It
is computationally stable because computing $\ell_n(\tau; \theta)$ for
one $\tau$ at a time is essentially one-dimensional problem.
Note that Carrasco \textit{et al.}~\cite{Caretal07} considered CCF-based
generalized method of moment estimation by considering a continuum of
$\tau$s in a functional space via covariance operator,
but the covariance operator may not be invertible due to zero eigenvalues.
Hence, Carrasco \textit{et al.}~\cite{Caretal07} needed ridging to avoid the
invertible issue, which makes the computation quite involved.

\section{Test for model specification}\label{sec3}

In this section we consider testing for the validity of (\ref
{eq:paralevy}) via testing for the parametric specification of the CCF
$\psi_t(u; \theta)$. Tests for model specification of a continuous-time
Markov process have been proposed by Chen, Gao and Tang~\cite{CheGaoTan08} and A\"
{\i}t-Sahalia,
Fan and Peng~\cite{AtSFanPen09}.
Despite the fact that parameter estimation based on the transition
density is asymptotically efficient, it is unclear
if a test based on the transition density is more powerful than one
based on the CCF.
The choice is clearer when the transition density does not admit a
closed form while the CCF does, since the latter is a test valid at any
level of the sampling interval $\delta$.

Let the underlying process that generates the observed sample path $\{
X_t\}_{t=1}^n$ be
\begin{equation}
\mathrm{d}X_t=\mu(X_t)\,\mathrm{d}t + \sigma(X_t) \,\mathrm{d} L_{t},
\label{eq:levy}
\end{equation}
whose CCF is $\psi(u; X_t)$. The process (\ref{eq:paralevy}) is a
parametric specification of (\ref{eq:levy}). To emphasize the
dependence of the CCF on $X_t$, we write in this section $\psi_t(u)$ as
$\psi(u, X_t)$, $\psi_t(u; \theta)$ as $\psi(u, X_t; \theta)$ and other
quantities in a similar fashion. We consider testing
\[
H_0\dvt P\{\psi_t(u) = \psi_t(u; \theta_0)\}=1 \qquad\mbox{for all $u \in
R^d$ and some $\theta_0 \in\Theta$},
\]
against a sequence of local alternative hypotheses
\[
H_1\dvt P\{\psi_t(u) = \psi_t(u; \theta_0) + c_n \Delta_n(u; X_t)\}=1
\qquad\mbox{for all $u \in R^d$},
\]
where $\{c_n\}$ is a sequence of non-random real constants converging
to zero at a certain rate, and $\{\Delta_n(u; X_t)\}$ is a sequence of
bounded complex functions which are continuous at $u =0$ and $\Delta
_n(0; X_t) \equiv0$; see Condition C6 in the \hyperref[app]{Appendix} for extra restrictions.

Since the target of inference is a conditional quantity, we need to
work with a kernel smoothed version of $\ell_n(\theta)$.
Let $K$ be a kernel function which is a symmetric probability density
in $R^d$, and~$h$ be a smoothing bandwidth that tends to $0$ as $n \to
\infty$.
A smoothed version of $L_n(\tau, \theta)$ is
\begin{equation}
L_{n h}(\tau, x; \theta) = \max\prod_{t=1}^n p_t(\tau, x) \label{eq:ELsm},
\end{equation}
subject to $\sum_{t=1}^n p_t(\tau, x ) = 1$ and $
\sum_{t=1}^n p_t(\tau,x) K_h(x-X_t) \vec{\varepsilon }(\tau, X_t;
\theta) = 0$.

Let $\ell_{n h} (\tau, x, \theta) = -2 \log\{ L_{n h}(\tau, x,
\theta)
n^{n} \}$ be the log-EL ratio.
Then the integrated log-EL ratio for $\theta$ is
\[
\ell_{n h}(\theta) =\int\int\ell_{n h} (\tau, x, \theta) \pi
_1(\tau)
\pi_2(x) \,\mathrm{d}\tau \,\mathrm{d}x ,
\]
where $\pi_1$ and $\pi_2$ are probability weight functions on the
frequency space and the state space, respectively. We can choose $\pi
_1$ to be the same as the $\pi$ in Section~\ref{sec3}.

The test statistic is $\ell_{n h}(\hat{\theta}_n)$, where
$\hat{\theta}_n$ is the empirical likelihood\ estimator proposed in
Section~\ref{sec3}. As a
matter of fact, we can employ any estimator with $n^{1/2}$-rate of convergence.
To appreciate the meaning of the test statistic,
let $W_h(x-X_t) = K_h(x-X_t) / \sum_{j=1}^n K_h(x-X_j)$ be the
Nadaraya--Watson kernel weight,
${\varepsilon }_{n,h}(\tau, x; \theta) = \sum_{t=1}^n K_h(x-X_t)
\varepsilon (\tau,
X_t; \theta)$ be the kernel smooth of the residuals, $\tilde
{{\varepsilon
}}_{n,h}(\tau, x; \theta) = ({\varepsilon }_{n h}(\tau, x; \theta),
{\varepsilon }_{n
h}(-\tau, x; \theta) )^T$ and $R(K) = \int K^2(t) \,\mathrm{d}t$. It can be shown
by a similar derivation in Chen, H\"ardle and Li~\cite{CheHarLi03} that
\begin{eqnarray}\label{eq:exp}
\ell_{nh}(\theta) &=& n h^d R^{-1}(K) \int\int\tilde{{\varepsilon }}_{n,
h}^{\star} (\tau, x; \theta) V^{-1}(\tau,x; \theta_0, \theta)
\tilde
{{\varepsilon }}_{n, h}(\tau, x; \theta)
\nonumber
\\[-10pt]
\\[-8pt]
\nonumber
& &\phantom{n h^d R^{-1}(K) \int\int}{} \times\pi_1(\tau) f(x)
\pi_2(x)
\,\mathrm{d}\tau \,\mathrm{d}x + \mathrm{O}_p\{(nh^d)^{-1/2} \log^3(n) + h^2 \log^2(n)\},
\end{eqnarray}
where $V(\tau,x;\theta_0, \theta) = \operatorname{Var}\{ \tilde{{\varepsilon }}
(\tau, X_t;
\theta) |X_t =x \}$, and $f(x)$ is the density of $X_t$. So, the test
statistic is asymptoticly equivalent to a $L_2$-measure of the averaged
``residuals'' $ \tilde{{\varepsilon }}_{n, h}^{\star} (\tau, x; \theta)$,
inversely weighted by the covariance matrix function $V$.
Hence the proposed test is similar in tune to Fan and Zhang~\cite{Fa03} for
testing diffusion processes, and of H\"ardle and Mammen~\cite{Ha93} and Wang
and Van Keilegom~\cite{Wa07} for testing regression functions.

We need the following notations to describe the power property. Let
$V(\tau_1, \tau_2, x) = E\{ \tilde{\varepsilon }(\tau_1, X_t;  \theta
_0) \tilde
{\varepsilon }^{\star}(\tau_2, X_t; \theta_0) | X_t =x\}$, then
$V(\tau, \tau,
x; \theta_0, \theta_0) = V(\tau, x)$, defined earlier.
Express the matrices
\[
V(\tau_1, \tau_2, x) = ( V_{lk}(\tau_1, \tau_2, x) )_{1 \le
l, k \le2} \qquad\mbox{and}\qquad V^{-1}(\tau, x) = (\nu^{lk}(\tau
,x) )_{1 \le l, k \le2}.
\]
Furthermore, we choose $c_n = n^{-1/2} h^{d/4}$ and define
\begin{eqnarray*}
\eta_n(\tau, X_t) &=& w(\tau;X_t) \Delta_n(u, X_t),\qquad \tilde{\eta
}_n(\tau, X_t) = (\eta_n(\tau, X_t), \eta_n(-\tau, X_t))^T,
\\
\mu_n &=& \int\int\tilde{\eta}_n^{\star} (\tau, x) V^{-1}(\tau, x;
\theta_0,\theta_0) \tilde{\eta}_n (\tau, x) \pi_1(\tau) \pi
_2(x) f(x) \,\mathrm{d} \tau \,\mathrm{d}x,
\end{eqnarray*}
$\sigma_n^2 = 2R^{-2}(K)h^{-d} \gamma^2(K, V, \pi_1, \pi_2)$ where
\begin{eqnarray}\label{eq:appengamma}
&&\gamma^2(K, V, \pi_1, \pi_2)\nonumber\\ 
&&\quad =
K^{(4)}(0) \int\int\int\sum_{l_1, k_1, l_2,k_2}^2 V_{l_1 l_2}(-\tau
_1, \tau_2, x) V_{k_1 k_2}(\tau_1, -\tau_2,x) \nu^{l_1, k_1}(\tau
_1, x)
 \\
&&\hspace*{20pt}\phantom{K^{(4)}(0) \int\int\int\sum_{l_1, k_1, l_2,k_2}^2}{}\times \nu^{l_2, k_2}(\tau_2, x) \pi
_1(\tau
_1) \pi_1(\tau_2) \pi_2^2(x) \,\mathrm{d}\tau_1 \,\mathrm{d} \tau_2 \,\mathrm{d} x,\nonumber
\end{eqnarray}
where $K^{(4)}$ is the $4$th convolution of the kernel function $K$.

The asymptotic normality of $\ell_{nh}(\hat{\theta}_n)$ is given in the
following theorem.

\begin{theorem}\label{th2} Under Conditions \textup{C1--C6} given in the \hyperref[app]{Appendix},
\begin{equation}
h^{-d/2} \bigl( \ell_{n h}(\hat{\theta}_n) - 2 - h^{d/2}
\mu_n \bigr)
\stackrel{d} \to N(0, 2R^{-2}(K)\gamma^2(K, V, \pi_1, \pi_2)).
\label{eq:clt}
\end{equation}
\end{theorem}

We note that $\mu_n=2$ under $H_0$.
Under $H_1$, since $\Delta_n(u, x)$ is non-vanishing with respect to
$u$, $\tilde{\eta}_n(\tau, x)$ is non-vanishing with respect to $u$ for
all $x$ in the support of $f$, which leads to a positive quantity $\mu
_n$, due to $V^{-1}(\tau, x; \theta_0,\theta_0)$ being a Hermitian matrix.
Since no restriction has been imposed on the functional form of $\Delta
_n(u, X_t)$, it means that the test is powerful for a wide range of
local alternatives. Indeed, if $\hat{\gamma}^2(K, {V}, \pi_1,\pi
_2)$ is
a consistent estimator of $\gamma^2(K, V, \pi_1, \pi_2)$,
the asymptotic normality-based test for $H_0$ with $\alpha$-level of
significance rejects $H_0$ if
\[
\ell_{n h}(\hat{\theta}_n) \ge2 + z_{1 - \alpha}\sqrt2 h^{d/2}
R^{-1}(K) \hat{\gamma}(K, {V}, \pi_1,\pi_2),
\]
where $z_{1 - \alpha}$ is the $1-\alpha$ quantile of the standard
normal distribution.
Theorem~\ref{th2} implies that the power of the test under $H_1$ is
\[
\Phi\biggl( - z_{1-\alpha} + { R(K) \mu_n \over\sqrt{ 2} \gamma(K, V,
\pi_1, \pi_2) } \biggr),
\]
where $\Phi$ is the standard normal distribution function.

It is known that the choice of bandwidth is important in any test based
on the kernel smoothing technique.
To make the test less sensitive to the choice of smoothing bandwidth,
we propose carrying out the test based on a set of bandwidths, say $\{
h_1, \ldots, h_k\}$, for a fixed integer $k$ such that $h_{i} = c_i h$
for some constants $c_1 < c_2 < \cdots < c_k$. Here $h$ is a reference
bandwidth which may be obtained via the cross-validation method.

This means that we have a set of the EL ratios
$\{\ell_{n h_1}(\hat{\theta}_n), \ldots, \ell_{n h_k}(\hat{\theta
}_n) \}
$ corresponding to the bandwidth set, and the overall test statistic is
\begin{equation}
T_n = \max_{1 \le i \le k} \bigl\{ h_i^{-d/2} \bigl(\ell_{n h_i}(\hat{\theta}_n)
-2 \bigr) \bigr\}. \label{eq:combined}
\end{equation}

To describe the asymptotic distribution of $T_n$, let $K^{(2)}(z, c) =
\int K(u) K(z + c u) \,\mathrm{d}u$ be a generalization to the convolution of $K$,
$\nu(t) = \int\{ K^{(2)}(t u,t) \}^2 \,\mathrm{d}u $ and
\[
\Sigma_J = {2 \over R^2(K)} \int\int\pi_1(\tau_1) \pi_1(\tau_2)
\pi
_2^2(x) \,\mathrm{d}x \,\mathrm{d} \tau_1 \,\mathrm{d} \tau_2 \bigl( (c_j/c_i)^{d} \nu(c_i/c_j)
\bigr)_{J\times J}.
\]

\begin{theorem}\label{th3} Under Conditions \textup{C1--C6}, $T_n \stackrel{d} \to
\max_{1
\le k \le J} Z_k$ as $ n \to\infty$,
where
\[
(Z_1, \dots, Z_J)^T \sim N(0, \Sigma_J).
\]
\end{theorem}

Let $t_{\alpha}$ be the $1-\alpha$ quantile of $T_{n}$, where $\alpha
\in(0,1)$ is the nominal size of the test.
The following parametric bootstrap procedure is employed to approximate
$t_{\alpha}$:

Step 1: Simulate a sample path $\{X_t^{\ast}\}_{t=1}^{n}$ at the same
frequency $\delta$ according to the model under~$H_0$ with the CCF
based estimate $\hat{\theta}_n$.

Step 2: Let $\tilde{\theta}_n^{\ast}$ be the estimate of $\theta$ under
$H_0$ using the resample path $\{X_t^{\ast}\}_{t=1}^{n}$ obtained in
Step~1, and $T_n^{\ast}$ be the version of $T_n$ for the resampled path.

Step 3: For a large positive integer $B$, repeat Steps 1 and 2 $B$ times
and obtain, after ranking,
$T_{n}^{(1) \ast} \le T_{n}^{(2) \ast} \le\cdots\le T_{n}^{(B) \ast}$.

Then, the Monte Carlo approximation of $t_{\alpha}$ is $T_{n}^{ ([B
(1-\alpha)] + 1) \ast}$. The proposed test rejects $H_0$ if
$T_{n}(\hat\theta_n) \ge T_{n}^{ ([B (1-\alpha)] + 1) \ast}$. The
justification of the above bootstrap procedure can be made based on
Theorem~\ref{th3} via the standard techniques for instance those given in Chen,
Gao and Tang~\cite{CheGaoTan08}.

\section{Simulation study}\label{sec4}

We report in this section the results from our simulation studies which
are designed to verify the proposed parameter estimator and model
testing procedure.
To evaluate the quality of the proposed EL estimator, we first chose
two univariate diffusion processes with known transition densities, so
that the MLEs can be compared with the proposed EL estimates. The two
processes are the Vasicek model (Vasicek~\cite{Vas77}) (VSK),
\begin{equation}\label{eq:VSK} \mathrm{d}X_{t} = \kappa(\alpha-X_{t})\,\mathrm{d}t +
\sigma
\,\mathrm{d}B_{t},
\end{equation}
and the Cox--Ingersoll--Ross model (Cox, Ingersoll and Ross~\cite{CoxIngRos85}) (CIR),
\begin{equation}\label{CIR}
\mathrm{d}X_{t} = \kappa(\alpha-X_{t})\,\mathrm{d}t + \sigma
\sqrt{X_{t}}\,\mathrm{d}B_{t},
\end{equation}
where $\kappa$, $\alpha$ and $\sigma$ are unknown parameters which
represent the mean reverting rate, long-run mean and volatility of the
process, respectively. Both processes are widely used in interest rate
modeling and various option price formulation. For the Vasicek model,
the transition distribution of $X_{t+1}|X_t$ is a normal distribution
$N(\alpha+(X_{t}-\alpha)\operatorname{exp}(-\kappa\delta), \sigma^2
(1-\operatorname{exp}(-2\kappa
\delta))/(2\kappa))$. For the CIR model, when $2\kappa\alpha/\sigma
^2 >
1$, $X_{t+1}|X_t$ is a multiple of a non-central Chi-square random
variable with degrees of freedom $4\kappa\alpha/\sigma^2$ and
non-centrality parameter $c X_{t}\operatorname{exp}(-\kappa\delta)$, where the
multiplier is $1/c$ with $c=4\kappa/(\sigma^2(1-\operatorname{exp}(-\kappa\delta)))$.
The CCFs of these two models can easily be derived from their known
transitional densities.

We then considered estimation for the jump diffusion model VSK-MJ as
given in (\ref{eq:VSK-MJ}) based on its CCF function
\begin{equation}
\psi_{t}(u;\theta)=\operatorname{exp}\biggl\{\frac{\sigma^2 u^2}{4\kappa}(\mathrm{e}^{-2\kappa
\delta
}-1)-\lambda\delta+\gamma+\mathrm{i}\bigl(\alpha u(1-\mathrm{e}^{-\kappa\delta})+u
\mathrm{e}^{-\kappa
\delta}X_{t}\bigr)\biggr\},
\label{eq:ccf_vsk_mj}
\end{equation}
where $\gamma=\lambda/(2\kappa) \int_{\mathrm{e}^{-2\kappa\delta}}^{1}
\operatorname{exp}(-\eta
^2 u^2 y/2)/y \,\mathrm{d}y$.
For comparison, we approximated its transition density by a mixture of
normal distributions, $(1-\lambda\delta)N(\mu_{\delta},\sigma
_{\delta
}^2)+\lambda\delta N(\mu_{\delta},\sigma_{\delta}^2+\eta^2),$
which is
a first order approximation proposed in A\"{\i}t-Sahalia, Fan and Peng
\cite{AtSFanPen09}. Here,
$\mu_{\delta}=\alpha+(X_{t}-\alpha)\operatorname{exp}(-\kappa\delta),$ and
$\sigma
_{\delta}^2=\sigma^2(1-\operatorname{exp}(-2\kappa\delta))/(2\kappa).$ The approximate
MLEs were obtained based on the mixture approximation given above.

We also consider the Inverse Gaussian OU process (IG-OU) in (\ref
{eq:IG-OU}), that is, the process $X_t$ follows the inverse Gaussian
law $\operatorname{IG}(a,b)$, for every $t$ when $X_0$ is generated from $\operatorname{IG}(a,b)$.
The CCF of this process is
\begin{equation}\psi_{t}(u;\theta)=\operatorname{exp}\bigl\{-a\bigl(\sqrt{-2\mathrm{i}u+b^2}-\sqrt
{-2\mathrm{i}u\mathrm{e}^{-\lambda
\delta}+b^2}\bigr)+\mathrm{i}u\mathrm{e}^{-\lambda\delta}X_{t}\bigr\}. \label{eq:IG-OUCCF}
\end{equation}
Since neither the exact transition density nor its approximation is
available, we were content with carrying out estimation with the
proposed methods.

The last simulation model considered for the estimation is a bivariate
extension of the univariate Ornstein--Uhlenbeck process (BI-OU),
\begin{equation} \label{eq:BI-OU} \mathrm{d}X_{t}={\bolds{\kappa}}({\bolds{\alpha}
}-X_{t})\,\mathrm{d}t+{\bolds{\sigma}} \,\mathrm{d}B_{t},
\end{equation}
where\vspace*{1pt} $X_{t}=(X_{1t},X_{2t})$,
$
{\bolds{\kappa}}=
({
{\kappa_{11}\enskip 0} \atop
{\kappa_{21} \enskip \kappa_{22}}}) ,
{\bolds{\alpha}}=
(
{
\alpha_{1} \atop
\alpha_{2}
}
)
\mbox{ and }
{\bolds{\sigma}}=
(
{
{\sigma_{11}\atop 0} \enskip
{0 \atop \sigma_{22}}}
).
$
Under the condition that the eigenvalues of the matrix ${\bolds{\kappa}}$
have positive real parts,
the process is stationary with transition distribution being a
bivariate normal $N(m(\delta,X_{t}),\Omega(\delta))$,
where $m(\delta,X_{t})=\alpha+\operatorname{exp}(-\kappa\delta)(X_{t}-\alpha)$,
$\Omega
(\delta)=\Sigma-\operatorname{exp}(-\kappa\delta) \Sigma \operatorname{exp}(-\kappa^{T} \delta
)$ and
\[
\Sigma=\frac{1}{2\operatorname{tr}(\kappa)\operatorname{Det}(\kappa)}\bigl\{\operatorname{Det}(\kappa)\sigma\sigma
^{T}+\lbrace\kappa-\operatorname{tr}(\kappa) \rbrace\sigma\sigma^{T} \lbrace
\kappa
-\operatorname{tr}(\kappa) \rbrace^{T}\bigr\}.
\]
The CCF of the process is known to be
$\psi_{t}(u_1,u_2; \theta)=\operatorname{exp}\lbrace \mathrm{i} u^{T} m(\delta,X_{t}) -u^T
\Omega(\delta) u/2 \rbrace$ for $u=(u_1,u_2)^T$.

We then carried out simulations to evaluate the ability of the proposed
tests in detecting model deviations. When we chose the simulation
models, we had in mind two issues in finance that have drawn
considerable research attention recently. The first issue is
whether the process is subject to jumps, and the second is whether we
could differentiate two processes with different jump rates. Our
simulation study formulated two settings of hypotheses to address these
two issues.
In the first setting, we tested
\[
\mbox{$H_0$}\dvt \mbox{The process is the VSK model}.
\]
In the second setting, we tested
\[
\mbox{$H_0$}\dvt \mbox{The process is the jump diffusion model VSK-MJ}.
\]
For computing the powers, in the first setting we used the data
simulated from $H_1\dvt$ the jump diffusion model VSK-MJ to test the null
model which does not have jumps; in the second setting, we used the
data simulated from $H_1\dvt$ the inverse Gaussian OU model which has
infinite-activity jumps to test the null hypothesis that prescribes a
finite-activity jump process.

For each model, we simulated 500 sample paths which were observed at
monthly observations ($\delta=1/12$) for $n=125, 250, 500$,
respectively. The choices of parameter values were motivated by Chen,
Gao and Tang~\cite{CheGaoTan08} and Ait-Sahalia, Fan and Peng~\cite{AtSFanPen09}.

In parameter estimation, we discovered that for both real and imaginary
parts of the CCF, their nonparametric smoothing estimators are
wave-like functions and roughly diminish to zero at the same points,
which creates a region denoted as $S_t$ (here the subscript $t$
indicates that the region depends on $X_t$). In practice, we searched
on a couple of grid points in the data range of~$X_t$ and picked the
union of $S_t$ as the support region $S$ for the frequency domain of
$\psi_t(u;\theta)$ in the estimation. We then chose the uniform density
as the weight function $\pi$ over the support region.

In model testing, a similar effort was initially made to obtain the
support region of the nonparametric CCF estimate, denoted as $S_{\mathrm{NP}}$,
and the support region of the theoretical CCF under $H_0$, denoted as
$S_{H_0}$. Here the theoretical CCF under $H_0$ used $\hat{\theta}_n$
from our EL method. Then the support region of the frequency domain in
testing was taken as the union of $S_{\mathrm{NP}}$ and $S_{H_{0}}$. We chose
the uniform density as the weight function over this support region for
testing. There is little contribution to the integrated empirical
likelihood ratio $\ell_{nh}(\hat{\theta}_n)$ from outside the support
region. The biweight Kernel $K(u)=15/16(1-u^2)^2 I(|u|\leq1)$ was used
for smoothing in testing. The bandwidth selection is described in
Section~\ref{sec3}. The bandwidth sets were specified in Tables~\ref{tb:test1}
and~\ref{tb:test2} for the two test settings. It is observed that the
values of the bandwidths were quite small, which was due to the rapid
oscillation of the CCF curves which favored smaller bandwidth in the
curve fitting.

We chose $w(u, r;X_t) = \mathrm{e}^{\mathrm{i} r^T X_t}$ throughout our simulation study
as it is the optimal instrument suggested in Carrasco \textit{et al.}
\cite{Caretal07}. Some numerical exploration (not reported) indicated the choice
of the function $w(\cdot)$ is not crucial in the context of the paper.
For testing, we picked the unit instrument to reduce computing burden.

%
\begin{table}
\caption{Empirical averages and their
standard errors (in parentheses) of the maximum (MLE) or approximate
maximum (AMLE) likelihood estimates and the proposed empirical
likelihood\ estimates
(EL) under the four univariate models}\label{tb:uni}
\begin{tabular*}{250pt}{@{\extracolsep{\fill}}lllll@{}}
\multicolumn{5}{c@{}}{(a) Vasicek model}\\
\hline
$n$& & $\kappa=0.858$ & $\alpha=0.089$ & $\sigma=0.047$\\
\hline
125 & MLE & 1.383 (0.603) &0.090 (0.015) &0.047 (0.003)\\
& EL & 1.305 (0.643) &0.090 (0.017) &0.046 (0.004)\\[3pt]
250 & MLE & 1.118 (0.397) &0.090 (0.011) &0.047 (0.002)\\
& EL & 1.052 (0.410) &0.089 (0.013) &0.046 (0.002)\\[3pt]
500 & MLE & 0.966 (0.240) &0.089 (0.008) &0.047 (0.002)\\
& EL & 0.951 (0.273) &0.089 (0.009) &0.047 (0.002)\\
\hline
\end{tabular*}\vspace*{9pt}
\begin{tabular*}{250pt}{@{\extracolsep{\fill}}lllll@{}}
\multicolumn{5}{c@{}}{(b) CIR model}\\
\hline
$n$& & $\kappa=0.892$ & $\alpha=0.091$ & $\sigma=0.181$\\
\hline
125 & MLE & 1.372 (0.644) &0.091 (0.019) &0.183 (0.012)\\
& EL & 1.290 (0.719) &0.093 (0.023) &0.178 (0.014)\\[3pt]
250 & MLE & 1.127 (0.374) &0.090 (0.013) &0.182 (0.008)\\
& EL & 1.089 (0.435) &0.091 (0.015) &0.179 (0.009)\\[3pt]
500 & MLE & 1.000 (0.245) &0.091 (0.010) &0.182 (0.006)\\
& EL & 0.977 (0.290) &0.092 (0.011) &0.180 (0.007)\\
\hline
\end{tabular*}\vspace*{9pt}
\begin{tabular*}{\textwidth}{@{\extracolsep{\fill}}lllllll@{}}
\multicolumn{7}{c@{}}{(c) Jump diffusion VSK-MJ model}\\
\hline
$n$& & $\kappa=0.858$ & $\alpha=0.089$ & $\sigma=0.047$ & $\lambda
=2.0$ &
$\eta=0.067$\\
\hline
125 & AMLE & 1.056 (0.381) &0.093 (0.020) &0.046 (0.005) & 1.770 (0.723) &
0.060 (0.016)\\
& EL & 1.090 (0.261) &0.084 (0.031) &0.048 (0.009) & 1.851 (0.323) &
0.066 (0.020)\\[3pt]
250 & AMLE & 0.977 (0.226) &0.093 (0.013) &0.047 (0.003) & 1.659 (0.466) &
0.059 (0.010)\\
& EL & 1.043 (0.201) &0.090 (0.023) &0.048 (0.007) & 1.825 (0.236) &
0.068 (0.015)\\[3pt]
500 & AMLE & 0.939 (0.145) &0.092 (0.009) &0.047 (0.002) & 1.620 (0.311) &
0.060 (0.007)\\
& EL & 1.018 (0.115) &0.089 (0.018) &0.049 (0.005) & 1.801 (0.163) &
0.068 (0.012)\\
\hline
\end{tabular*}\vspace*{9pt}
\begin{tabular*}{250pt}{@{\extracolsep{\fill}}lllll@{}}
\multicolumn{5}{c@{}}{(d) Inverse Gaussian OU model}\\
\hline
$n$& & $\lambda=10.0$ & $a=1.0$ & $b=20.0$\\
\hline
125 & EL & 10.328 (3.665) &1.048 (0.106) &20.722 (2.146)\\
250 & EL & 11.154 (1.976) &1.059 (0.043) &21.380 (0.878)\\
500 & EL & 11.489 (1.652) &1.031 (0.024) &20.846 (0.461)\\
\hline
\end{tabular*}
\end{table}

Table~\ref{tb:uni} reports the empirical averages of the parameter
estimates and their standard errors as well as the true parameter
values used for simulation. When the sample size increases, standard
errors of all the proposed estimates decrease, indicating the
consistency of the estimators.
We observe from Table~\ref{tb:uni}(a)--(b) for the VSK and CIR models
where the MLEs are available, the proposed EL estimates are quite close
to the MLEs. Although the EL estimates tend to have larger standard
errors than the MLEs, we do note that under the VSK model in Table \ref
{tb:uni}(a), the bias of EL estimates for the mean reverting parameter
$\kappa$ are smaller than the corresponding MLEs for all $n=125$,
$n=250$ and $n=500$.
For the jump diffusion model VSK-MJ (Table~\ref{tb:uni}(c)), we see
the EL estimates are consistently more efficient than the approximate
MLEs in the estimation of $\kappa$ and the Poisson intensity $\lambda$.
For the Inverse Gaussian OU model, which does not have the MLE to
compare with, the proposed estimates as reported in Table~\ref{tb:uni}(d) are
close to the true values, and the standard errors converge as the
sample size increases.

\begin{table}
\caption{Empirical averages and their
standard errors (in parentheses) of the maximum (MLE) likelihood
estimates and the proposed empirical likelihood\ estimates (EL) under
the Bivariate OU model}\label{tb:bi}
\begin{tabular*}{\textwidth}{@{\extracolsep{\fill}}lllll@{}}
\hline
$n$& & $\kappa_{11}=0.22$ & $\kappa_{21}=0.2$ & $\kappa_{22}=0.5$\\
\hline
125 & MLE & 0.441 (0.197) &0.395 (0.270) &0.607 (0.176)\\
& EL & 0.381 (0.208) &0.525 (0.238) &0.594 (0.192)\\[3pt]
250 & MLE & 0.353 (0.165) &0.307 (0.148) &0.563 (0.110)\\
& EL & 0.354 (0.178) &0.449 (0.184) &0.564 (0.153)\\[3pt]
500 & MLE & 0.280 (0.118) &0.241 (0.104) &0.526 (0.068)\\
& EL & 0.261 (0.168) &0.383 (0.154) &0.487 (0.112)\\
\hline
\end{tabular*}\vspace*{9pt}
\begin{tabular*}{\textwidth}{@{\extracolsep{\fill}}llllll@{}}
\hline
$n$& & $\alpha_{1}=0.08$ & $\alpha_2=0.09$ & $\sigma_{11}=0.09$&
$\sigma
_{22}=0.17$\\
\hline
125 & MLE & 0.145 (0.166) &0.099 (0.056) &0.167 (0.067) &0.080 (0.079)\\
& EL & 0.141 (0.141) &0.117 (0.085) &0.129 (0.044) &0.071 (0.034)\\[3pt]
250 & MLE & 0.141 (0.151) &0.096 (0.036) &0.140 (0.065) &0.116 (0.074)\\
& EL & 0.142 (0.129) &0.094 (0.073) &0.095 (0.033) &0.094 (0.028)\\ [3pt]
500 & MLE & 0.102 (0.120) &0.092 (0.023) &0.115 (0.051) &0.146 (0.055)\\
& EL & 0.099 (0.108) &0.104 (0.064) &0.077 (0.024) &0.105 (0.028)\\
\hline
\end{tabular*}
\end{table}

Table~\ref{tb:bi} reports the estimates for the bivariate OU process
and shows that the EL estimates are close to the corresponding MLEs,
providing the further evidence of the effectiveness of our EL estimator
for multivariate process estimation. We also found that the EL
estimates for the long run mean $\alpha_1$ and the volatility $\sigma
_{11}$ of the first process have smaller biases and standard errors
than the MLEs for all $n=125$, $n=250$ and $n=500$.

%
\begin{table}
\caption{$H_0$: VSK versus $H_1$: the
jump diffusion model VSK-MJ}\label{tb:test1}
\begin{tabular*}{\textwidth}{@{\extracolsep{\fill}}lld{1.3}d{1.3}d{1.3}d{1.3}d{1.3}l@{}}
\multicolumn{8}{c@{}}{(a) Size evaluation (in percentage)}\\
\hline
$n=125$&Bandwidth&0.012&0.017&0.021&0.025&0.030&\multicolumn{1}{l@{}}{Overall}\\
&Size&4.6&5.6&5.4&5.8&5.6&4.8\\[3pt]
$n=250$&Bandwidth&0.012&0.015&0.018&0.021&0.024&\multicolumn{1}{l@{}}{Overall}\\
&Size&5.6&6.2&6.2&6.0&5.8&5.4\\[3pt]
$n=500$&Bandwidth&0.011&0.013&0.015&0.018&0.020&\multicolumn{1}{l@{}}{Overall}\\
&Size&5.0&5.6&5.6&5.4&5.6&5.0\\
\hline
\end{tabular*}\vspace*{9pt}
\begin{tabular*}{\textwidth}{@{\extracolsep{\fill}}lld{2.3}d{2.3}d{2.3}d{2.3}d{2.3}l@{}}
\multicolumn{8}{c@{}}{(b) Power evaluation (in percentage)}\\
\hline
$n=125$&Bandwidth&0.016&0.021&0.026&0.032&0.037&\multicolumn{1}{l@{}}{Overall}\\
&Power&72.0&71.6&70.4&69.2&65.8&72.2\\[3pt]
$n=250$&Bandwidth&0.016&0.019&0.022&0.026&0.029&\multicolumn{1}{l@{}}{Overall}\\
&Power&82.4&82.4&82.2&82.4&82.2&82.6\\[3pt]
$n=500$&Bandwidth&0.014&0.017&0.019&0.021&0.024&\multicolumn{1}{l@{}}{Overall}\\
&Power& 95.0&94.8&94.6&94.4&94.2&94.8\\
\hline
\end{tabular*}
\end{table}
%

%
%
\begin{table}
\caption{$H_0$: the jump diffusion model
VSK-MJ versus $H_1\dvt$ the inverse Gaussian OU model}\label{tb:test2}
\begin{tabular*}{\textwidth}{@{\extracolsep{\fill}}lld{1.3}d{1.3}d{1.3}d{1.3}d{1.3}l@{}}
\multicolumn{8}{c@{}}{(a) Size evaluation (in percentage)}\\
\hline
$n=125$&Bandwidth&0.017&0.022&0.028&0.034&0.040&\multicolumn{1}{l@{}}{Overall}\\
&Size&3.4&3.6&4.0&3.6&4.6&4.6\\[3pt]
$n=250$&Bandwidth&0.017&0.021&0.024&0.028&0.032&\multicolumn{1}{l@{}}{Overall}\\
&Size&4.6&4.6&4.6&4.6&5.0&4.8\\[3pt]
$n=500$&Bandwidth&0.016&0.019&0.021&0.024&0.026&\multicolumn{1}{l@{}}{Overall}\\
&Size&5.0&5.2&5.2&5.0&5.0&5.0\\
\hline
\end{tabular*}\vspace*{9pt}
\begin{tabular*}{\textwidth}{@{\extracolsep{\fill}}lld{2.3}d{2.3}d{2.3}d{2.3}d{2.3}l@{}}
\multicolumn{8}{c@{}}{(b) Power evaluation (in percentage)}
\\
\hline
$n=125$&Bandwidth&0.008&0.012&0.017&0.021&0.026&\multicolumn{1}{l@{}}{Overall}\\
&Power&71.6&73.8&73.2&71.4&71.2&74.4\\[3pt]
$n=250$&Bandwidth&0.008&0.011&0.014&0.017&0.020&\multicolumn{1}{l@{}}{Overall}\\
&Power&84.0&84.2&83.4&81.8&81.4&84.4\\[3pt]
$n=500$&Bandwidth&0.008&0.010&0.012&0.014&0.016&\multicolumn{1}{l@{}}{Overall}\\
&Power& 90.1&88.9&89.5&85.1&85.4&90.2\\
\hline
\end{tabular*}
\end{table}

Tables~\ref{tb:test1} and~\ref{tb:test2} report the empirical size and
power of the proposed test based on $B=250$ bootstrap resampled paths
for each simulation. They contain the sizes and powers for the overall
test that is based on the five bandwidth set, and for the tests that
only use one bandwidth. We observe that the tests gave satisfactory
sizes under both testing settings. In the first test where we used the
data from the jump diffusion model VSK-MJ to test the continuous
diffusion model VSK, the powers range from $65\%$ to $95\%$ across the
different sample sizes and bandwidths. In the second test where we used
data simulated from the infinity-activity jump process (the inverse
Gaussian OU) to test the finite-activity jump process (the jump
diffusion VSK-MJ), the powers range from $71\%$ to $90\%$ across the
different sample sizes and bandwidth choices.

We also compared our methods with Carrasco \textit{et al.}~\cite{Caretal07} for
estimation, and with Chen, Gao and Tang~\cite{CheGaoTan08} for testing. To save
space, we reported the results in details in the supplemental article
(Chen, Peng and Yu~\cite{ChePenYu}).

\section{A case study}\label{sec5}

In this section, we examine empirically the capability of our testing
procedure in detecting jumps using the secondary market quotes of the
3-month Treasury Bill (T-bill) between January 1, 1965 and February 2,
1999. This bill was sampled at monthly frequency, and in total we had
410 observations.
The mean of these bills is $0.065$, the volatility is $0.026$, the mean
of the differences is very close to zero ($1.5\times10^{-5}$) and the
standard deviation of the differences is $0.005$.
The sample period contains some large movements that turn out to
coincide with arrivals of macroeconomic news (Johannes~\cite{Joh04}). The goal
of this empirical study was to test whether the underlying process is
subject to jumps or not.

\begin{table}
\caption{Empirical estimation for the
3-month T-bill Data}\label{tb:real_est}
\begin{tabular*}{200pt}{@{\extracolsep{\fill}}ld{1.5}d{1.5}d{1.5}@{}}
\multicolumn{4}{c@{}}{(a) VSK model}\\
\hline
&\multicolumn{1}{l}{\hspace*{-3pt}$\kappa$}&\multicolumn{1}{l}{\hspace*{-3pt}$\alpha$}&\multicolumn{1}{l@{}}{\hspace*{-3pt}$\sigma$}\\
\hline
MLE&0.277&0.065&0.019\\
&(0.1800) & (0.0117) & (0.0007) \\[3pt]
EL&0.274&0.059&0.018\\
& (0.1956) & (0.0136) &(0.0007) \\
\hline
\end{tabular*}\vspace*{9pt}
\begin{tabular*}{200pt}{@{\extracolsep{\fill}}ld{1.5}d{1.5}d{1.5}@{}}
\multicolumn{4}{c@{}}{(b) CIR model}
\\
\hline
&\multicolumn{1}{l}{\hspace*{-3pt}$\kappa$}&\multicolumn{1}{l}{\hspace*{-3pt}$\alpha$}&\multicolumn{1}{l@{}}{\hspace*{-3pt}$\sigma$}\\
\hline
MLE&0.182&0.066&0.061\\
&(0.1697) &(0.0179) &(0.0021)\\[3pt]
EL&0.182&0.064&0.057\\
& (0.1934) &(0.0374)&(0.0021) \\
\hline
\end{tabular*}\vspace*{9pt}
\begin{tabular*}{\textwidth}{@{\extracolsep{\fill}}ld{1.5}d{1.5}d{1.5}d{1.5}d{1.5}@{}}
\multicolumn{6}{c@{}}{(c) VSK-MJ model}\\
\hline
&\multicolumn{1}{l}{\hspace*{-3pt}$\kappa$}&\multicolumn{1}{l}{\hspace*{-3pt}$\alpha$}&\multicolumn{1}{l}{\hspace*{-3pt}$\sigma$}&
\multicolumn{1}{l}{\hspace*{-3pt}$\lambda$}&\multicolumn{1}{l@{}}{\hspace*{-3pt}$\eta$}\\
\hline
AMLE&0.071&0.077&0.009&1.863&0.012\\
&(0.0170)& (0.0129)&(0.0004)&(0.3282)&(0.0015)\\[3pt]
EL&0.072&0.076&0.008&1.862&0.013\\
&(0.0143)& (0.0136)&(0.0008)& (0.1569)&(0.0021) \\
\hline
\end{tabular*}\vspace*{9pt}
\begin{tabular*}{200pt}{@{\extracolsep{\fill}}ld{1.5}d{1.5}d{1.5}@{}}
\multicolumn{4}{c@{}}{(d) Inverse Gaussian OU model}
\\
\hline
&\multicolumn{1}{l}{\hspace*{-3pt}$\lambda$}&\multicolumn{1}{l}{\hspace*{-3pt}$a$}&\multicolumn{1}{l@{}}{\hspace*{-3pt}$b$}\\
\hline
EL&0.264&1.139&12.558\\
&(0.0342)&(0.1364)&(0.8970)\\
\hline
\end{tabular*}
\end{table}

The proposed parameter estimates under each of the four univariate
models considered in the simulation study are reported in Table \ref
{tb:real_est}. For comparison, the MLEs or the approximate MLEs are
also reported except for the Inverse Gaussian OU model. For the
univariate diffusion models VSK and CIR, and the jump diffusion model
VSK-MJ, the proposed parameter estimates based on CCF are very similar
to the MLEs or the approximate MLEs. The EL estimates of the long-run
mean $\alpha$ are $0.059$ for VSK and $0.064$ for CIR, both of which
are close to the summary statistic of mean rates ($0.065$). In VSK, the
average volatility of 3-month T-bill monthly return (difference) is
estimated to be $\sigma\sqrt{\delta}=0.018 \sqrt{1/12}=0.005$, which
is also close to the summary statistic of volatility for the change
($0.005$). However the conditional volatility of monthly change in CIR
model is $\sigma\sqrt{\delta X_{t}}$, and $X_t$ has a long-run
average $0.064$ which is
less than $1$. Therefore, the process needs to have higher $\sigma$
($0.057$) to bring up the average volatility of monthly change to the
same level reflected by the real data. In the jump diffusion model
VSK-MJ, our estimate of $\lambda$ suggests on average about 2 jumps per
year. Relative to VSK and CIR models, the estimate for parameter
$\sigma
$ in the jump diffusion VSK-MJ model is much smaller ($0.008$),
indicating that allowing jumps in the process helps to capture large
movements in the interest rate, and, as a result, the continuous part
of the process does not have to be as volatile as the one in VSK or CIR models.

%
\begin{table}
\caption{p-values for the 3-month
T-bill data}\label{tb:real_test}
\begin{tabular*}{\textwidth}{@{\extracolsep{\fill}}lld{2.3}d{2.3}d{2.3}d{2.3}d{2.3}d{2.4}@{}}
\hline
&& \multicolumn{6}{l@{}}{Bandwidth}\\[-6pt]
&& \multicolumn{6}{c@{}}{\hrulefill}\\
&&\multicolumn{1}{l}{0.010}&\multicolumn{1}{l}{0.012}&
\multicolumn{1}{l}{0.014}&\multicolumn{1}{l}{0.016}&
\multicolumn{1}{l}{0.018}&\multicolumn{1}{l@{}}{Overall}\\
\hline
VSK&Test Stats& 21.971&19.225&16.145&13.267&10.786&14.828\\
&$l_{0.05}^{*}$&3.228&3.123&2.845&2.724&2.647&1.462\\
&p-values&0.0&0.0&0.0&0.0&0.0&0.0\\[3pt]
CIR&Test Stats&6.015&4.775&3.755&2.954&2.335&3.546\\
&$l_{0.05}^{*}$&2.782&2.739&2.825&2.650&2.448&1.229\\
&p-values&0.0&0.01&0.02&0.026&0.054&0.0\\[3pt]
VSK-MJ&Test Stats&37.204&40.901&45.046&49.878&55.561&25.600\\
&$l_{0.05}^{*}$&35.669&43.548&52.247&62.744&74.298&28.751\\
&p-values&0.046&0.074&0.102&0.126&0.148&0.0880\\[3pt]
IG-OU &Test Stats&10.716&9.374&7.962&6.663&5.528&6.870\\
&$l_{0.05}^{*}$&40.463&47.665&46.444&42.396&41.750&27.940\\
&p-values&0.11&0.148&0.124&0.128&0.122&0.162\\
\hline
\end{tabular*}
\end{table}

We then applied the proposed test for the validity of each of the four
models. The bandwidth prescribed by the CV was $0.01$. By exploring the
kernel estimators of the CCF, a reasonable range for $h$ was from
$0.01$ to $0.018$, that offered smoothness from slightly
under-smoothing to slightly over-smoothing. The bandwidth range used in
our empirical study consisted of five equally spaced bandwidths ranging
from $0.01$ to $0.018$. Table~\ref{tb:real_test} reports p-values of
single bandwidth and the overall tests for the four models. There is no
empirical support for the VSK model. The CIR model performs a little
bit better as the distances between the test statistics and the
critical values decrease, but the model is still rejected at
significance level of $0.05$ in the overall test and almost all the
single bandwidth tests. We can not reject the jump diffusion model
VSK-MJ in the overall test and the single bandwidth tests except the
one with the smallest bandwidth (p-value $=0.046$). This constitutes a
strong indication of the presence of jumps and implies that adding
(finite-activity) jumps does help to capture the underlying dynamics of
the interest rates. By allowing the infinite-activity jumps in the
models, the p-values of the tests for the inverse Gaussian OU model are
very supportive, even for the small bandwidths, suggesting that the
infinite-activity jump model might potentially model the dynamics of
the 3-month T-bill rates better. A possible reason for this is that the
jump diffusion model VSK-MJ can only generate small continuous
movements from Brownian motion and big spikes from the compound Poisson
component, but it could miss the movements that are between (i.e., the
movements with median sizes). However, the inverse Gaussian OU process
is more flexible since it can generate small, median and big movements
with infinite arrival rates; therefore it could fill in a gap in the
VSK-MJ model by capturing movements that are too large for Brownian
motion to model but too small for the compound Poisson process to capture.

\begin{appendix}
\section*{Appendix}\label{app}

The following conditions are required in our analysis.

\begin{longlist}
\item[C1.] {The stochastic processes given in (\ref{eq:paralevy}) and (\ref
{eq:levy}) admit unique weak solution respectively, which are $\alpha$-mixing with mixing coefficient
$\alpha(t)=C \mathrm{e}^{-\lambda t}$ where $\alpha(t)=\sup\{|P(A\cap B) -
P(A)P(B)|\dvt A\in\Omega_1^s, B\in\Omega_{s+t}^{\infty}\}$
for all $s,t\geq1$, where $C$ is a finite positive constant, and
$\Omega_i^j$ denotes the $\sigma$-field generated by $\{X_t\dvt i\leq
t\leq j\}$}.

\item[C2 (Smoothness).]  $\psi_t(\tau; \theta)=: \psi(\tau; \theta, X_t)$ and
$E\{ \varepsilon_t(\tau; \theta)\}$ are third continuous differentiable
with respect to $\theta$ within a neighborhood of $\theta_0$ which is
defined in C3. $\pi(\cdot)$ is a bounded probability density supported
on a compact set $S \subset R^d$; and the diffusion function $\sigma
(x)$ is positive definite.

\item[C3.] The parameter space $\Theta$ is an open subset of $R^p$, and the
true parameter $\theta_0$ is the unique root of $E\{\varepsilon_t(\tau;
\theta)\}=0$ for all $\tau\in S$; and for any $\theta_1 \ne\theta_2$,
$P\{ \psi_t(\cdot; \theta_1) \ne\psi_t(\cdot; \theta_2, X_t)
\} > 0$.

\item[C4 (Invertibility).] The Hermitian matrix $\operatorname{Var}\{ \tilde{\varepsilon
}_t(\tau
; \theta_0) \}$ is positive definite almost everywhere for $\tau\in
R^{2d}$ with respect to the Lebesgue measure in $R^{2d}$;
$\Gamma(\theta_0)$ defined in (\ref{eq:Gamm1}) is invertible.

\item[C5.] The kernel $K(\cdot)$ is a $r$th order symmetric kernel supported
on $[-1,1]^d$ and has bounded second derivative. We assume $d < 4$ and
the smoothing bandwidth $h = \mathrm{O}\{n^{-1/(d+2r)}\}$.
The bandwidth set $\{h_1, \ldots, h_k\}$ satisfies $h_i = c_i h$ for
constants $c_i$ such that $c_1 < c_2 < \cdots < c_k$ where $k$ is an
integer not depending on $n$.

\item[C6.] $\{\Delta_n(u; X_t)\}$ is a sequence of complex functions
continuous at $u =0$ and $\Delta_n(0; X_t) \equiv0$, $\sup_{n}
|\Delta
_n(u; X_t)| \le M_1$ almost surely and the Lebesgue measure of $\{ u |
\Delta_n(u, x) \ne0 \}$ is positive for all $x$ in the support of the
marginal density $f$, and $c_n =n^{-1/2} h^{-d/4}$ which is the order
of the difference between $H_0$ and $H_1$.
\end{longlist}

We need C1 as the basic condition for the stochastic processes
involved. Ait-Sahalia~\cite{AtS96} and
Genon-Catalot, Jeantheau and Laredo~\cite{GenJeaLar00} provide conditions on the
underlying processes such that
Assumption C1 held. In particular, Ait-Sahalia~\cite{AtS96} provides
conditions so that the observed sequences
are $\beta$-mixing, which is automatically $\alpha$-mixing. We require
that the rate of decay is exponentially fast
to simplify the technical arguments.
C2 consists of smoothness conditions regarding the CCFs and C3 is for
identification of parameters. C4 ensures the covariance matrix is
invertible, which is easier to be justified for our low-dimensional
formulation of estimation and testing approaches. C5 on the kernel and
bandwidth are standard in nonparametric curve estimation. The
assumption of $d < 4$ is to make the bias in the kernel estimation a
smaller order of $h^{d/2}$ so that the bias is stochastically
negligible relative to $\ell_{n h}(\theta_0)$. The kernel method will
encounter the curse of dimensionality when $d \ge4$. Also, the
commonly used processes in finance and other stochastic modeling tend
to have dimension less than 4.
The bandwidth selected by either cross validation or the plug-in method
satisfies the order specified in C5.
The first part of C6 regarding $\Delta_n(u; X_t)$ is to qualify $\psi
_t(u; \theta)$ under $H_1$ as a bona fide characteristic function,
whereas the part that requires positive measure on the set $\{ u |
\Delta_n(u, x) \ne0 \}$ is to make $H_1$ a genuine sequence of
alternative hypotheses.

\begin{pf*}{Proof of Lemma \protect\ref{le1}}
By combining results in Kitamura~\cite{Kit97} and Chen, H\"ardle and Li
\cite{CheHarLi03} for the empirical likelihood of $\alpha$-mixing processes, we
can show that
\begin{equation}
\label{proof1-1}
\lambda(\tau;\theta)=A_n^{-1}(\tau;\theta)\Biggl\{\frac1n\sum
_{t=1}^n\vec
{\varepsilon }_t(\tau;\theta)\Biggr\}+\mathrm{o}(n^{-1/3})=\mathrm{O}(n^{-1/3})
\end{equation}
almost surely and
uniformly in $\Vert \theta-\theta_0\Vert \le n^{-1/3}$ and $\tau^T\in S$.
Denote $\theta=\theta_0+un^{-1/3}$. It follows from (\ref{proof1-1})
and Taylor's expansion that, uniformly in $\Vert u\Vert =1$,
\begin{eqnarray}\label{proof1-2}
&&\ell_n(\theta)\nonumber\\
&&\quad= \int\Biggl\{2\sum_{t=1}^n\lambda^T(\tau;\theta)\vec{\varepsilon
}_t(\tau;\theta)
-\sum_{t=1}^n\{\lambda^T(\tau;\theta)\vec{\varepsilon }_t(\tau
;\theta)\}^2\Biggr\}\pi
(\tau) \,\mathrm{d}\tau+\mathrm{o}(n^{1/3})\nonumber \\
&&\quad=\int n\Biggl\{\frac1n\sum_{t=1}^n\vec{\varepsilon }^T_t(\tau;\theta
_0)+\frac1n\sum
_{t=1}^n\frac{\partial\vec{\varepsilon }^T_t(\tau;\theta
_0)}{\partial\theta
}un^{-1/3}\Biggr\} A_n^{-1}(\tau;\theta) \\
&&{}\qquad\phantom{\int}{} \times\Biggl\{\frac1n\sum_{t=1}^n\vec{\varepsilon }^T_t(\tau;\theta
_0)+\frac
1n\sum_{t=1}^n\frac{\partial\vec{\varepsilon }^T_t(\tau;\theta
_0)}{\partial
\theta}un^{-1/3}\Biggr\}\pi(\tau) \,\mathrm{d}\tau+\mathrm{o}(n^{1/3})\nonumber\\
&&\quad=  \int n\biggl\{E\biggl(\frac{\partial\vec{\varepsilon }^T_1(\tau;\theta
_0)}{\partial
\theta}\biggr)un^{-1/3}\bigl(1+\mathrm{o}(1)\bigr)\biggr\}A^{-1}(\tau,\tau;\theta_0,\theta_0) \nonumber\\
&&\qquad\phantom{\int}{} \times\biggl\{E\biggl(\frac{\partial\vec{\varepsilon }_1(\tau;\theta
_0)}{\partial\theta
}\biggr)un^{-1/3}\bigl(1+\mathrm{o}(1)\bigr)\biggr\}\pi(\tau) \,\mathrm{d}\tau+\mathrm{o}(n^{1/3})\nonumber\\
&&\quad\ge \frac12cn^{1/3}\nonumber
\end{eqnarray}
almost surely,
where $c>0$ is the smallest eigenvalue of
\[
\sup_{\tau\in S}E\Biggl(\frac{\partial\vec{\varepsilon }^T_1(\tau;\theta
_0)}{\partial
\theta}\Biggr)A^{-1}(\tau, \tau;\theta_0,\theta_0) E\Biggl(\frac{\partial
\vec{\varepsilon
}_1(\tau;\theta_0)}{\partial\theta}\Biggr).
\]
Similarly,
\begin{eqnarray*}
\ell_n(\theta_0)&=&\int\Biggl\{\sum_{t=1}^n\vec{\varepsilon }^T_t(\tau
;\theta_0)\Biggr\}
A^{-1}(\tau, \tau;\theta_0,\theta_0)\Biggl\{\frac1n\sum_{t=1}^n\vec
{\varepsilon
}_t(\tau;\theta_0)\Biggr\}\pi(\tau) \,\mathrm{d}\tau+\mathrm{o}(1)\\
&=&\mathrm{o}(n^{1/3}),
\nonumber
\end{eqnarray*}
almost surely.
This, together with (\ref{proof1-2}), implies that $\ell_n(\theta)$ has
a minimum value in the interior of the ball $\Vert \theta-\theta_0\Vert \le
n^{-1/3}$, and this value satisfies $\frac{\partial}{\partial\theta}
\ell_n(\theta)=0$, that is, the second equation in (\ref{th1-1}) by
noting (\ref{eq:3.5.2}). The first equation follows directly from
(\ref{eq:3.5.2}).
\end{pf*}

\begin{pf*}{Proof of Theorem \protect\ref{th1}} It follows from limit theorems for
martingale difference that
\begin{equation}
\label{proof2-1}
\cases{
\displaystyle\frac{\partial}{\partial\theta}Q_{1n}(\tau;\theta_0,0)=\frac
1n\sum
_{t=1}^n\frac{\partial}{\partial\theta}\vec{\varepsilon }_t(\tau
;\theta
_0)\stackrel{p}{\to}M_0E\biggl\{\frac{\partial}{\partial\theta}\tilde
\varepsilon
_1(\tau;\theta_0)\biggr\},\vspace*{2pt}\cr
\displaystyle\frac{\partial}{\partial\lambda^T}Q_{1n}(\tau;\theta_0,0)=-\frac
1n\sum
_{t=1}^n\vec{\varepsilon }_t(\tau;\theta_0)\vec{\varepsilon }^T_t(\tau
;\theta
_0)\stackrel{p}{\to} -M_0A(\tau,\tau;\theta_0,\theta_0)M_0^{\star
},\vspace*{2pt}\cr
\displaystyle\frac{\partial}{\partial\theta}Q_{2n}(\tau;\theta_0,0)=0,\vspace*{2pt}\cr
\displaystyle\frac{\partial}{\partial\lambda^T}Q_{2n}(\tau;\theta_0,0)=\frac
1n\sum
_{t=1}^n\frac{\partial}{\partial\theta}\vec{\varepsilon }^T_t(\tau
;\theta
_0)\stackrel{p}{\to}E\biggl\{\frac{\partial}{\partial\theta}\tilde
\varepsilon
^{\star}_1(\tau;\theta_0)\biggr\}M_0^{\star}}\
\end{equation}
uniformly in $\tau^T\in S$.
Put $\delta_n=\Vert \hat\theta_n-\theta_0\Vert +\sup_{\tau^T\in S}
\Vert \lambda(\tau;\hat\theta_n)\Vert $. Then it follows from Taylor's
expansion that
\begin{eqnarray}
\label{proof2-2}
0&=&Q_{1n}(\tau;\hat\theta_n,\lambda(\tau;\hat\theta_n))
\nonumber
\\[-8pt]
\\[-8pt]
\nonumber
&=&Q_{1n}(\tau;\theta_0,0)+\frac{\partial Q_{1n}(\tau;\theta
_0,0)}{\partial\theta}(\hat\theta_n-\theta_0)+\frac{\partial
Q_{1n}(\tau
;\theta_0,0)}{\partial\lambda^T}\lambda(\tau;\hat\theta
_n)+\mathrm{o}_p(\delta_n)
\end{eqnarray}
uniformly in $\tau^T\in S$, and
\begin{eqnarray}\label{proof2-3}
0&=&\int Q_{2n}(\tau;\hat\theta_n,\lambda(\tau
;\hat
\theta_n))\pi(\tau) \,\mathrm{d}\tau \nonumber\\
&=&\int\biggl\{Q_{2n}(\tau;\theta_0,0)+\frac{\partial Q_{2n}(\tau;\theta
_0,0)}{\partial\theta}(\hat\theta_n-\theta_0)
+\frac{\partial Q_{2n}(\tau;\theta_0,0)}{\partial\lambda^T}\lambda
(\tau
;\hat\theta_n)\biggr\}\pi(\tau) \,\mathrm{d}\tau\qquad \\
& & {}+\mathrm{o}_p(\delta_n).\nonumber
\end{eqnarray}
By (\ref{proof2-1})--(\ref{proof2-3}), we have
\begin{eqnarray} \label{proof2-4}
&& \hat\theta_n-\theta_0
\nonumber
\\[-8pt]
\\[-8pt]
\nonumber
&&\quad=-\Gamma^{-1}(\theta_0)\int E\biggl\{\frac{\partial}{\partial\theta
}\tilde
\varepsilon^{\star}_1(\tau;\theta_0)\biggr\}A^{-1}(\tau;\theta_0,\theta
_0)M_0^{-1}\frac1n\sum_{t=1}^n\vec{\varepsilon }_{t}(\tau;\theta
_0)\pi(\tau)
\,\mathrm{d}\tau+\mathrm{o}_p(\delta_n).\qquad
\end{eqnarray}
Hence the theorem follows from (\ref{proof2-4}) and the central limit
theorem for Martingale difference.
\end{pf*}

\begin{pf*}{Proof of Theorem \protect\ref{th2}}
Define $V(\tau_1, \tau_2, x; \theta_0, \theta) = E\{ \tilde
{\varepsilon }(\tau
_1, X_t; \theta) \tilde{\varepsilon }^{\star}(\tau_2, X_t; \theta) |
X_t =x\}$
and write $V(\tau, x; \theta_0, \theta) = V(\tau, \tau, x; \theta_0,
\theta)$.
Since $\hat{\theta}_n$ is $\sqrt{n}$-consistent to $\theta_0$, we have
\begin{eqnarray} \label{eq:exptest}
\ell_{n h}(\hat{\theta}_n) &=& \ell_{n h, 1}(\theta_0) + n h^d
R^{-1}(K) \lbrace(\hat{\theta} - \theta_0)^T S_{n, h}(\theta_0) +
S_{n, h}^{\star}(\theta_0) (\hat{\theta}_n - \theta_0) \nonumber\\
& &\hspace*{100pt}{} + (\hat{\theta}_n - \theta_0)^T \Gamma_{n, h}(\theta_0) (\hat
{\theta
}_n - \theta_0) \rbrace\\
&&{}+ \mathrm{O}_p\{(nh^d)^{-1/2} \log^3(n)+ h^2
\log^2(n)\},\nonumber
\end{eqnarray}
where
\begin{eqnarray}\label{eq:Gammanh}
\ell_{n h, 1}(\theta_0) &=& n h^d R^{-1}(K) \int\int\tilde
{{\varepsilon
}}_{n, h}^{\star} (\tau, X_t; \theta_0) V^{-1}(\tau,x; \theta
_0,\theta
_0)\nonumber\\[-1pt]
&&\phantom{n h^d R^{-1}(K) \int\int}{} \times\tilde{{\varepsilon }}_{n, h}(\tau, x; \theta_0) \pi_1(\tau)
f^{-1}(x) \pi_2(x) \,\mathrm{d}\tau \,\mathrm{d}x,
\nonumber\\[-1pt]
S_{n,h}(\theta_0) &=& \int\int{\partial\tilde{{\varepsilon
}}_{n,h}^{\star
}(\tau, x; \theta_0) \over\partial\theta} V^{-1}(\tau,x; \theta
_0,\theta_0) \tilde{{\varepsilon }}_{n,h}(\tau, x; \theta_0)\nonumber\\[-1pt]
& &\phantom{\int\int}{} \times\pi_1(\tau) \pi_2(x) f^{-1}(x) \,\mathrm{d} \tau \,\mathrm{d}x,
\nonumber\\[-1pt]
\Gamma_{n h}(\theta_0) &=& \int\int{\partial\tilde{\varepsilon
}^{\star
}_{n,h}(\tau, x; {\theta}_0) \over\partial\theta} V^{-1}(\tau,x;
\theta_0,\theta_0) {\partial\tilde{\varepsilon }_{n, h}(\tau, x;
{\theta}_{0})
\over\partial\theta}
\nonumber
\\[-9pt]
\\[-9pt]
\nonumber
&&\phantom{\int\int}{} \times\pi_1(\tau) \pi_2(x) f^{-1}(x) \,\mathrm{d} \tau \,\mathrm{d}x.
\end{eqnarray}
As $S_{n,h}(\theta_0) = \mathrm{O}_p(n^{-1/2})$,
\begin{eqnarray}
\ell_{n h}(\hat{\theta}_n) &=& \ell_{n h, 1}(\theta_0) + \mathrm{O}_p\{
(nh^d)^{-1/2} \log^3(n) + h^2 \log^2(n) + h^d\}. \label{eq:Aellnh1}
\end{eqnarray}
Note that
\begin{eqnarray}\label{eq:A.8}
& & \ell_{n h, 1}({\theta}_0)\nonumber\\[-1pt]
&&\quad = nh^d R^{-1}(K) \int\int n^{-1} \sum_{t_1=1}^n K_h(x-X_{t_1})
\lbrace\tilde{\varepsilon }^{\star}(\tau, X_{t_1}) + c_n \tilde{\eta
}_n^{\star
} (\tau, X_{t_1}) \rbrace\nonumber\\[-1pt]
&&\qquad\phantom{nh^d R^{-1}(K) \int\int}{} \times V^{-1}(\tau, x;\theta_0,\theta_0) n^{-1} \sum_{t_2=1}^n
K_h(x-X_{t_2})
\nonumber
\\[-9pt]
\\[-9pt]
\nonumber
&&\qquad\phantom{nh^d R^{-1}(K) \int\int}{} \times \{\tilde{\varepsilon } (\tau, X_{t_2}) + c_n
\tilde{\eta
}_n(\tau, X_{t_2}) \}\\[-1pt]
&&\qquad\phantom{nh^d R^{-1}(K) \int\int} {}\times\pi_1(\tau) \pi_2(x) f^{-1}(x) \,\mathrm{d} \tau \,\mathrm{d}x +\mathrm{o}_p(h^{d/2})\nonumber \\[-1pt]
&&\quad= R^{-1}(K)( H_{n1} + H_{n2} + H_{n3} + H_{n4} ) +
\mathrm{o}_p(h^{d/2}),\nonumber
\end{eqnarray}
where, with the choice of $c_n=n^{-1/2} h^{-d/4}$,
\begin{eqnarray*}
H_{n 1} &=& n^{-1} h^d \sum_{t_1 \ne t_2} \int\int K_h(x-X_{t_1})
K_h(x-X_{t_2}) \tilde{\varepsilon }^{\star}(\tau, X_{t_1}) V^{-1}(\tau
, x)\\[-1pt]
&& \hspace*{66pt}{}\times\tilde{\varepsilon } (\tau, X_{t_2}) \pi_1(\tau) \pi_2(x)
f^{-1}(x) \,\mathrm{d} \tau \,\mathrm{d}x, \nonumber \\[-1pt]
H_{n 2} &=& n^{-1} h^d \sum_{t=1}^n \int\int K_h^2(x-X_{t}) \tilde
{\varepsilon
}^{\star}(\tau, X_{t}) V^{-1}(\tau, x) \tilde{\varepsilon } (\tau,
X_{t})\\[-1pt]
&& \hspace*{66pt}{}\times\pi_1(\tau) \pi_2(x) f^{-1}(x) \,\mathrm{d} \tau \,\mathrm{d}x, \\[-1pt]
H_{n 3} &=& 2n^{1/2} h^{3d/4} \int\int\tilde{\eta}_n^{\star} (\tau,
x) V^{-1}(\tau, x) n^{-1} \sum_{t=1}^n K_h(x-X_{t}) \tilde{\varepsilon
} (\tau,
X_{t})\\[-1pt]
&&\hspace*{66pt}{} \times\pi_1(\tau) \pi_2(x)f^{-1}(x) \,\mathrm{d} \tau \,\mathrm{d}x, \\[-1pt]
H_{n 4} &=& h^{d/2} \int\int\tilde{\eta}_n^{\star} (\tau, x)
V^{-1}(\tau, x) \tilde{\eta}_n (\tau, x) \pi_1(\tau) \pi_2(x) f^{-1}(x)
\,\mathrm{d} \tau \,\mathrm{d}x.
\end{eqnarray*}

We note that $H_{n2} = 2 R(K) + \mathrm{o}_p(h^d)$
and the integral in $H_{n3}$ is $\mathrm{O}_p(n^{-1/2})$. Hence, $H_{n 3} =
\mathrm{O}_p(n^{3d/4}) = \mathrm{o}_p(h^{d/2})$.

Now consider $H_{n 1}$. Clearly, $E(H_{n 1}) =0$ and the double
summation in $H_{n 1}$ constitutes a generalized $U$-statistic of order
two with the kernel
\begin{eqnarray*}
\xi_{t_1, t_2} &=&
\int\int K_h(x-X_{t_1}) K_h(x-X_{t_2}) \tilde{\varepsilon }^{\star
}(\tau,
X_{t_1}) V^{-1}(\tau, x;\theta_0,\theta_0) \tilde{\varepsilon } (\tau
, X_{t_2})
\\[-1pt]
&&\hspace*{22pt}{}\times\pi_1(\tau) \pi_2(x) f^{-1}(x) \,\mathrm{d} \tau \,\mathrm{d}x.
\end{eqnarray*}
The $U$-statistic is degenerate, due to $\{ \tilde{\varepsilon } (\tau,
X_{t_2}) \}$ being martingale differences.

Let $\sigma_n^2 = \sum_{ 1 \le t_1 \ne t_2 \le n} \sigma_{t_1, t_2}^2$
where $\sigma_{t_1, t_2}^2 = \operatorname{Var}(\xi_{t_1, t_2})$.
Then, applying the central limit theorem for generalized $U$-statistics
for $\alpha$-mixing sequences (Gao and King~\cite{GaoKin}),
we have
\begin{equation}
\sigma_n^{-1} \sum_{t_1 \ne t_2} \xi_{t_1, t_2} \stackrel{d} \to
N(0,1). \label{eq:cltU}
\end{equation}
Furthermore, it can be shown, for instance, by following the route of
Chen, Gao and Tang~\cite{CheGaoTan08},
that $\sigma_n^2 = 2 n^2 \sigma_{n 0}^2 \{1+ \mathrm{o}(1)\}$ where
$\sigma_{n 0}^2 = E_{t_1} E_{t_2} (\xi_{t_1, t_2}^2)$. Here $E_{t_i}$
denote marginal expectation with respect to $(X_{t_i} , X_{t_i +1})$.

It can be shown that
\begin{eqnarray}\label{eq:appenvar1}
\sigma_{n 0}^2
&=&\int\!\!\int\!\!\int\!\!\int E_{t_1} E_{t_2} \Biggl\{ K_h(x_1-X_{t_1})
K_h(x_1-X_{t_2}) K_h(x_2-X_{t_1})K_h(x_2-X_{t_2}) \nonumber\\[-1pt]
&&\hspace*{76pt}{} \times  \sum_{l_1, k_1, l_2,k_2}^2 {\varepsilon
}_{l_1}(\tau_1, X_{t_1}) \varepsilon _{k_1}(\tau_1,X_{t_2}) {\varepsilon
}_{l_2}(\tau
_2, X_{t_1})\nonumber\\[-1pt]
&&\hspace*{121pt}{} \times\varepsilon _{k_2}(\tau_2,X_{t_2}) \nu^{l_1, k_1}(\tau_1,
x_1) \nu
^{l_2, k_2}(\tau_2, x_2)\Biggr\}\nonumber\\[-1pt]
&&\hspace*{42pt}{} \times\pi_1(\tau_1) \pi_1(\tau_2) f^{-1}(x_1) f^{-1} (x_2) \pi
_2(x_1) \pi_2(x_2) \,\mathrm{d}\tau_1 \,\mathrm{d} \tau_2 \,\mathrm{d} x_1\, \mathrm{d} x_2\nonumber\\[-1pt]
&=& \int\!\!\int\!\!\int\!\!\int E_{t_1} E_{t_2}\Biggl\{ K_h(x_1-X_{t_1})
K_h(x_1-X_{t_2}) K_h(x_2-X_{t_1})K_h(x_2-X_{t_2})\qquad\\[-1pt]
&&\hspace*{76pt}{} \times  \sum_{l_1, k_1, l_2,k_2}^2 V_{l_1
l_2}(-\tau
_1, \tau_2, X_{t_1}) V_{k_1 k_2}(\tau_1, -\tau_2, X_{t_2})\nonumber\\[-1pt]
&&\hspace*{121pt}{} \times\nu^{l_1, k_1}(\tau_1, x_1) \nu^{l_2, k_2}(\tau_2, x_2)
\Biggr\}\nonumber\\[-1pt]
&&\hspace*{42pt}{} \times\pi_1(\tau_1) \pi_1(\tau_2) f^{-1}(x_1) f^{-1}(x_2)\pi_2(x_1)
\pi_2(x_2) \,\mathrm{d}\tau_1 \,\mathrm{d} \tau_2 \,\mathrm{d} x_1 \,\mathrm{d} x_2\nonumber\\[-1pt]
&=& h^{-d} \gamma^2(K, V, \pi_1, \pi_2) \{ 1 +
\mathrm{O}(h^2)\},\nonumber
\end{eqnarray}
where $\gamma^2(K, V, \pi_1, \pi_2)$ is defined in (\ref{eq:appengamma}).
From (\ref{eq:cltU}) and (\ref{eq:appenvar1}), we have
\begin{equation}
h^{-d/2} H_{n1} \stackrel{d} \to N(0, 2\gamma^2(K, V, \pi_1, \pi_2) ).
\end{equation}
This, together with the results on $H_{n2}$ and $H_{n3}$, leads to
\begin{equation}h^{-d/2} \bigl( \ell_{n h}(\hat{\theta}) - 2 - \mu_n \bigr)
\stackrel{d} \to
N(0, 2R^{-2}(K)\gamma^2(K, V, \pi_1, \pi_2)) \label{eq:clt1},
\end{equation}
where $\mu_n = H_{n4}$. This completes the proof of Theorem~\ref{th2}.
\end{pf*}

\begin{pf*}{Proof of Theorem \protect\ref{th3}}
The proof can be made by applying the Cram\'er--Wold device and the
same technique in the proof of Theorem 2, followed by the mapping theorem.
\end{pf*}
\end{appendix}

\section*{Acknowledgements}
The authors thank two reviewers and an associate editor for helpful
comments, and acknowledge support from National Science Foundation
Grants DMS-06-04563, DMS-05-18904, SES-0631608 and DMS-10-05336. The first
author was partially supported through the Center for Statistical
Science at Peking University.

\begin{supplement} [id=suppA]
\stitle{Comparisons in estimation and testing with other methods\\}
\slink[doi]{10.3150/11-BEJ400SUPP}
\sdatatype{.pdf}
\sfilename{bej400\_supp.pdf}
\sdescription{We compared our methods with Carrasco \textit{et al.} \cite{Caretal07}
for estimation, and with Chen, Gao and Tang~\cite{CheGaoTan08} for testing. The
supplemental article (Chen, Peng and Yu~\cite{ChePenYu}) provides additional
tables from these comparisons.}
\end{supplement}


\printhistory

\end{document}